# An adaptive edge-based smoothed finite element method (ES-FEM) for phase-field modeling of fractures at large deformations

Fucheng Tian[a,b], Xiaoliang Tang[a,b], Tingyu Xu[a,b], Liangbin Li[a,b,c*]

a *National Synchrotron Radiation Lab, University of Science and Technology of China, Hefei 230029, China*

b *Anhui Provincial Engineering Laboratory of Advanced Functional Polymer Film, University of Science and Technology of China, Hefei 230029, China* c *CAS Key Laboratory of Soft Matter Chemistry, University of Science and Technology of China, Hefei 230029, China*

**Abstract**

This work presents the Griffith-type phase-field formation at large deformation in the framework of adaptive edge-based smoothed finite element method (ES-FEM) for the first time. Therein the phase-field modeling of fractures has attracted widespread interest by virtue of its outstanding performance in dealing with complex cracks. The ES-FEM is an excellent member of the S-FEM family developed in combination with meshless ideas and finite element method (FEM), which is characterized by higher accuracy, 'softer' stiffness, and insensitive to mesh distortion. Given that, the advantages of the phase-field method (PFM) and ES-FEM are fully combined by the approach proposed in this paper. With the costly computational overhead of PFM and ES-FEM in mind, a well-designed multi-level adaptive mesh strategy was developed, which considerably improved the computational efficiency. Furthermore, the detailed numerical implementation for the coupling of PFM and ES-FEM is outlined. Several representative numerical examples were recalculated based on the proposed method, and its effectiveness is verified by comparison with the results in experiments and literature. In particular, an experiment in which cracks deflected in rubber due to impinging on a weak interface was firstly reproduced.

*Keywords*: ES-FEM; Adaptive; Phase-field method; Hyperelastic; large deformation

* Correspondence author: lbli@ustc.edu.cn, Tel & Fax: +86-551-63602081

# 1. Introduction

Hyperelastic materials like rubber, hydrogels, etc. have been widely used in industrial as well as academic research on the strength of their appealing properties such as high stretchability and reversibility [1-4]. In light of this, the prediction of crack initiation, propagation, and failure of such materials is of great significance for engineering applications. Early studies mostly utilized the discontinuous approaches to simulate the evolution of cracks [5-9]. However, this is still an intractable task for nonlinear elastomers that typically undergo complicated deformations. Along with the development of a category of the diffusive crack models [10-12], this predicament is expected to be solved. In contrast with the discontinuous crack modeling, such an approach does not require explicit tracking of sharp crack surfaces, which facilitates the handling of complex crack patterns as branching and intersecting [13, 14]. The phase field method (PFM), as one of the most prominent of these approaches, has spawned extensive research and is utilized in this investigation [15-23]. Relying on this model, the crack path is automatically determined by the principle of total potential energy minimization without ad hoc assumptions.

The phase field approach of fracture has two self-contained geneses, namely Griffith's theory and Ginzburg-Landau phase transition theory [24-28]. We consider a Griffith-type phase field model that originated in the mechanics community. This model roots in the variational approach of brittle fracture, a milestone work of Francfort and Marigo [28]. To set forth a numerically solvable version, Bourdin et al. introduced Ambrosio-Tortorelli regularizations into the variational form, which yielded the pristine phase field models of fracture [10, 27]. Thereafter, numerous researchers committed to perfecting the phase field model, leading to fruitful achievements [20-23, 29-31]. For instance, in the framework of thermodynamic consistency, Miehe et al. reconstructed the phase field model in a more stable form and proposed a spectral decomposition for strain tensor to avoid the nonphysical crack growth caused by compression [23, 30]. The vast majority of the previous reports are limited to small deformations. Hitherto, however, research on phase field modeling of fractures at large deformations is still deficient. To the best of our knowledge, Miehe et al. were the first to investigate the phase field modeling of fractures in rubber-like polymers at finite deformations [32]. Subsequently, Ambati et al., Borden et al. and Miehe et al. extended this model to ductile fracture [21, 33, 34]. Still based on Miehe's model, Loew et al. further introduced a rate-dependent phase-field damage model and demonstrated experimental verification [35]. More recently, Tang et al. proposed a novel strain energy-based decomposition scheme for general nonlinear elastic materials to distinguish the contribution of tension and compression to crack nucleation and propagation [36, 37].

Noted that the previous research related to the phase field approach to fracture is almost entirely conducted in the system of FEM. In spite of FEM is already one of the most popular methodologies for solving various partial differential equation (PDE) problems, its key thought emerged in the middle of the last century [38, 39]. Some innate drawbacks lead to the standard FEM is not the optimal choice at large deformations. For these reasons, Liu et al. developed a category of smooth finite element methods (S-FEM) which is an elaborate combination of standard FEM and some meshfree techniques [40-43]. The core of S-FEM is the strain smoothing technique, which stems from the stabilized conforming nodal integration (SCNI) concept proposed by Chen et al. [44]. Currently, S-FEM has evolved into a large family. Classified as per the construction of the smooth domain, which mainly comprises cell-based S-FEM (CS-FEM) [43], node-based S-FEM (NS-FEM) [45], edge-based S-FEM (ES-FEM) [46], and face-based S-FEM (FS-FEM) [47]. These different types of S-FEMs have been proven to hold their own unique properties. In general, the stiffness matrix produced by S-FEM is softer than FEM, which attenuates the overestimation of stiffness in FEM, resulting in higher accuracy and convergence rate [41]. Since no mapping operation is performed in S-FEM, the Jacobian matrix that is sensitive to mesh distortion does not exist. S-FEM is consequently more robust in handling mesh distortion and extreme deformation [40, 48, 49]. Besides, S-FEM can be constructed directly based on existing finite element meshes without adding additional degrees of freedom. With its distinctive attributes, S-FEM has been extensively applied in various aspects of solid mechanics, especially in terms of large deformation [48-53]. The latest report has reformatted the phase field model of fracture based on CS-FEM, nevertheless, it is still confined to the small deformation [54].

Apart from the aforesaid advantages, the primary deficiency of S-FEM lies in the high computational cost required, on account of its larger bandwidth of stiffness matrix than FEM [41]. Phase field approach also requires enormously fine meshes to accurately identify the crack path [13, 55, 56]. Conceivably, the combination of PFM and S-FEM is a daunting challenge for computing resources. In reality, multifarious sophisticated adaptive grid schemes have been developed to meet the exacting demands of PFM for mesh, e.g. hybrid adaptive mesh and multiscale mesh [55, 57]. However, most of the existing algorithms are designed for the phase field modeling of fracture in the framework of FEM [55, 56, 58]. To effectively perform the phase field modeling of the fracture at large deformations in the framework of S-FEM, developing an adaptive mesh scheme for the coupling of S-FEM and PFM is essential.

In this research, we present the formation of the phase field approach for modeling the fracture at large deformation in the framework of ES-FEM, one of the most promising S-FEM [46]. Thanks to the high accuracy and mesh distortion insensitivity of ES-FEM, the advantages of the phase field approach in dealing with complex cracks are sufficiently released in terms of large deformation. In consequence of the high computational cost of PFM and ES-FEM, a multi-level adaptive mesh scheme ameliorated from the *ha*-PFM we proposed earlier was further developed [57]. Distinguish with the adaptive grid in the context of FEM, adaptive ES-FEM requires to identify the connectivity of the edges and their supporting elements after each update of the mesh. Notwithstanding this strategy adds about 1%~2% extra computational overhead, it not only improves accuracy but also raises computational efficiency by roughly 20 times. On top of that, the outstanding performance of the presented method is testified by four representative benchmarks derived from experiments and previous simulations.

The rest of this manuscript is organized as follows. In Section 2, we briefly introduce the fundamental theory of finite (large) deformation and the hyperelastic Neo-Hookean constitutive model. On this basis, the governing equations of the phase field modeling of fracture at large deformation are derived. Section 3 outlines the theoretical aspects of ES-FEM containing the crucial strain smoothing technique and a novel adaptive grid algorithm. Detailed numerical implementations involving the discretization and linearization of weak form are derived in Section 4. Several validating examples are demonstrated in Section 5. At last, Section 6 summarizes this paper with some conclusions.

**2. Phase-field formation of fracture at large deformations**

In this section, a Griffith-type phase field model is formatted in the framework of large deformations. For clarity, we organize this section as the following four subsections. In subsection 2.1, starting with the basic kinematics, the concept of large deformation was introduced first. A hyperelastic Neo-Hookean model employed herein is described in subsection 2.2. Afterward, a thermodynamically consistent phase-field model of diffuse cracks rooted in fracture mechanics is outlined in subsection 2.3. Performing a variation on the Lagrangian yields the final governing equations, which is presented in subsection 2.4.

2.1. Kinematics

In the large deformation context, the initial configuration and the current configuration require a clear distinction. We consider an arbitrary elastomer with an initial (undeformed) configuration of $\Omega_0$, in which the position vector of a material point $P_i$ is represented by **X**. During the deformation process, the motion of **X** is identified by mapping $\mathbf{x} = \chi(\mathbf{X}, t)$ in the current (deformed) configuration $\Omega$. Thus, the fundamental deformation gradient tensor can be expressed as [38]

$$\mathbf{F} = \nabla_{\mathbf{X}}\chi(\mathbf{X},t) = \mathbf{I} + \nabla_{\mathbf{X}}\mathbf{u}, \tag{1}$$

where **I** is the second-order unity tensor and **u** denotes the displacement field. The determinant of **F**, i.e. Jacobian $J = \det(\mathbf{F}) > 0$ establishes an integral mapping between the initial configuration ($\Omega_0$) and the current configuration ($\Omega$), which is written as

$$d\Omega = J d\Omega_0. \tag{2}$$

In the theoretic frames of nonlinear continuum mechanics, an essential deformation measure is the Green-Lagrange strain tensor, that is

$$\mathbf{E} = \frac{1}{2}(\mathbf{C} - \mathbf{I}). \tag{3}$$

Here, **C** is the right Cauchy-Green tensor in terms of the material coordinates, it follows that

$$\mathbf{C} = \mathbf{F}^T \mathbf{F} = (\nabla_{\mathbf{X}}\mathbf{u})^T + \nabla_{\mathbf{X}}\mathbf{u} + \nabla_{\mathbf{X}}\mathbf{u} \cdot (\nabla_{\mathbf{X}}\mathbf{u})^T + \mathbf{I}. \tag{4}$$

This deformation measurement is widely utilized in constitutive equations, like the hyperelastic Neo-Hookean model presented in the next subsection.

2.2. Hyperelastic model

Abundant hyperelastic constitutive models have been developed and commendably modeled the mechanical response of materials like rubber, hydrogel, etc. [38, 59, 60]. In this contribution, we consider a category of the isotropic elastomers characterized by the Neo-Hookean model. The free energy density in the absence of damage can be written as

$$\psi_0(\mathbf{F}) = \frac{\mu}{2}[\text{tr}[\mathbf{F}^T\mathbf{F}] - 3] + \frac{\mu}{\beta}[(J^{-\beta} - 1)] \quad (5)$$

with $\beta = 2\nu/(1-\nu)$ [32]. Where $\mu$ is the shear modulus, and $\nu$ stands for the Poisson ratio. Accordingly, we have the first Piola-Kirchhoff stress (PK1)

$$\mathbf{P}_0 = \partial_{\mathbf{F}}\psi_0(\mathbf{F}) = \mu[\mathbf{F} - J^{-\beta}\mathbf{F}^{-T}] \quad (6)$$

and the second Piola-Kirchhoff stress (PK2)

$$\mathbf{S}_0 = \mathbf{F}^{-1} \cdot \mathbf{P}_0 = \mu(\mathbf{I} - \mathbf{J}^{-\beta}\mathbf{C}^{-1}). \quad (7)$$

To facilitate the application of the Voight notation in programming, the tangent modulus is derived based on Eq. 7 as

$$\mathbb{C}^{SE} = 2\frac{\partial \mathbf{S}}{\partial \mathbf{C}} = \mu J^{-\beta}[C_{il}^{-1}C_{kj}^{-1} + C_{ik}^{-1}C_{jl}^{-1}] + \beta\mu J^{-\beta}\, C_{ij}^{-1}C_{kl}^{-1}. \quad (8)$$

Note that the derivation of the aforementioned formats is based on non-damaging materials. After introducing the phase-field damage variable, all the above formulas require to be multiplied by a monotonically increasing degradation function that will be presented later.

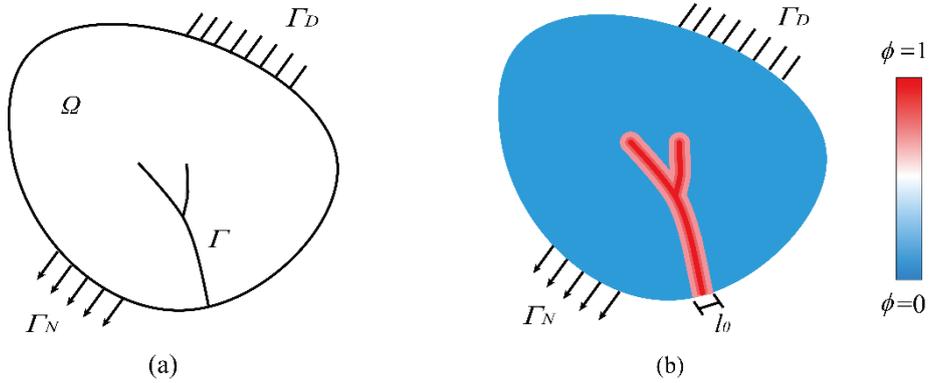

FIG. 1. (a) An elastomer containing a sharp crack $\Gamma$. (b) Phase-field approximate of diffuse cracks.

2.3. Phase-field representation of diffuse cracks

The quintessence of phase-field idea is to regularize the sharp crack topology by a limited diffuse damage band. Let us consider an arbitrary elastomer $\Omega \subset \mathbb{R}^n, (n=2,3)$ containing a sharp crack $\Gamma$ as the reference configuration and its surface is denoted by $\partial \Omega \subset \mathbb{R}^{n-1}$, as depicted in FIG. 1(a). To regularize the discontinuous crack surface, a time-dependent phase-field damage variable $\phi(\mathbf{X},t) \in [0,1]$ is introduced, where $\phi=0$ indicates undamaged and $\phi=1$ stands for the cracks (see FIG. 1(b)). According to the variational phase-field approach, the crack surface energy $\Psi(\Gamma)$ can be approximated by

$$\Psi(\Gamma) = \int_\Gamma G_c \mathrm{d}A \approx \int_\Omega G_c \gamma(\phi, \nabla\phi) \mathrm{d}A, \tag{9}$$

where $G_c$ is the fracture energy of materials and $\gamma$ is the crack surface density function. Note that the format of $\gamma$ is not unique. Throughout this work, two representative formats of $\gamma$ denoted by AT$_1$

$$\gamma(\phi, \nabla\phi) = \frac{3}{8}[\frac{1}{l_0}\phi + l_0(\nabla\phi \cdot \nabla\phi)] \tag{10}$$

and AT$_2$

$$\gamma(\phi, \nabla\phi) = \frac{1}{2}[\frac{1}{l_0}\phi^2 + l_0(\nabla\phi \cdot \nabla\phi)], \tag{11}$$

respectively, are considered [13, 61]. Herein, $l_0$ represents the regularization parameter, relating to the material's characteristic length $l_{ch}$. The discrepancy between AT$_1$ and AT$_2$ can be referred to the work of Tanné et al [61].

Previous experiments have indicated that a small portion of the potential energy is dissipated at the crack tip in the form of acoustic emission or heat generation [62-64].

For this, an extra dissipative term is added to the total potential energy functional $\Pi$, which is reformatted as

$$\underbrace{\Pi(\mathbf{u},\phi)}_{total\ potential} = \underbrace{\int_\Omega g(\phi)\psi_0(\mathbf{F})\mathrm{d}V}_{strain\ energy} + \underbrace{\int_\Omega G_c\gamma(\phi,\nabla\phi)\mathrm{d}V}_{crack\ surfacce\ energy} + \underbrace{\frac{\eta}{2}\int_\Omega(\frac{\partial\phi}{\partial t})^2\mathrm{d}V}_{dissipation} \\ -\underbrace{(\int_\Omega \bar{\mathbf{b}}\cdot\mathbf{u}\mathrm{d}V + \int_{\partial\Omega}\bar{\mathbf{t}}\cdot\mathbf{u}\mathrm{d}A)}_{external}, \quad (12)$$

where, $g(\phi)=(1-\phi)^2+k$ is a degradation function describing the attenuation of stored elastic energy due to crack evolution [23, 30]. $\bar{\mathbf{b}}$ and $\bar{\mathbf{t}}$ are the volume and traction force vectors, respectively. Furthermore, a positive parameter $k$ ($0<k<<1$) is added for numerical reasons.

2.4. Governing equations

By performing variational operations on the total energy functional $\Pi$, the strong form of governing equations for the phase-field description of the fracture at large deformations can be derived as

$$\begin{cases} \mathrm{Div}[\mathbf{P}]+\bar{\mathbf{b}}=\mathbf{0} \\ \dfrac{G_c}{l_0}[\phi-l_0^2\Delta\phi]+\eta\dot{\phi}=2(1-\phi)\psi_0 \end{cases} \quad (13)$$

with the boundary conditions

$$\begin{cases} \mathbf{P}\cdot\mathbf{n}=\bar{\mathbf{t}} & \text{on } \Gamma_N \\ \nabla\phi\cdot\mathbf{n}=0 & \text{on } \partial\Omega \end{cases}. \quad (14)$$

Note that the viscosity regularization is introduced into the above governing equations, which yields a rate-dependent crack growth model [32, 35]. Evidently, the general rate-independence can be effortlessly restored via setting $\eta=0$.

With respect to the Neo-Hookean model, the first Piola-Kirchhoff stress involving the phase-field damage reads

$$\mathbf{P} = ((1-\phi)^2 + k)\mu[\mathbf{F} - J^{-\beta}\mathbf{F}^{-T}]. \tag{15}$$

In general, cracks cannot heal (excluding certain hydrogels [65, 66]), for which we enforce $\frac{\partial \phi}{\partial t} \geq 0$ to meet the irreversibility of the crack. Besides, this restriction results in a nonnegative dissipation term

$$\mathcal{D} = \frac{\eta}{2} \int_{\Omega} (\frac{\partial \phi}{\partial t})^2 d\Omega \geq 0, \tag{16}$$

which explicitly satisfies the thermodynamic consistency condition [30, 67]. It is worth mentioning that the introduction of the dissipation term can stabilize the numerical calculations.

## 3. Theoretical aspects of ES-FEM

In this section, we briefly introduce the basic concept of S-FEM, i.e. strain (gradient) smoothing technique [40]. According to the dissimilar types of smoothing domains, several models with different features are proposed in the S-FEM family. Among these models, a very prominent ES-FEM is utilized to discretize the governing equations of phase field modeling for fracture at large deformations.

3.1. Strain smoothing

In the standard FEM, the entire solution domain $\Omega$ is divided into a set of elements $\Omega_e$. The compatible strain field is calculated by the gradient of the displacement field at the element level. And the global stiffness matrix can be obtained based on the assembly of the elements. However, in the S-FEM, the solution domain $\Omega$ is further split into a set of non-overlapped smoothing domains $\Omega_k^s (k=1,2,...N_{sd})$, where $N_{sd}$ is the total number of smoothing domains. Accordingly, the smoothed strain at the material point $\mathbf{X}_{sd}^s$ in $\Omega_k^s$ is defined by [41]

$$\bar{\boldsymbol{\varepsilon}}_k(\mathbf{X}_{sd}^s) = \int_{\Omega^h} \boldsymbol{\varepsilon}^h \Phi_k(\mathbf{X}-\mathbf{X}_{sd}^s) \mathrm{d}\Omega, \tag{17}$$

where $\boldsymbol{\varepsilon}^h$ is the compatible strain tensor, $\Phi_k(\mathbf{X}-\mathbf{X}_{sd}^s)$ is a smoothing function that satisfies

$$\Phi_k(\mathbf{X}-\mathbf{X}_{sd}^s) \text{ and } \int_{\Omega_k^s} \Phi_k(\mathbf{X}-\mathbf{X}_{sd}^s) \mathrm{d}\Omega = 1. \tag{18}$$

For the sake of simplicity, a Heaviside-type function is adopted [40], which is given by

$$\Phi_k(\mathbf{X}-\mathbf{X}_{sd}^s) = \begin{cases} 1/A_k, \mathbf{X} \in \Omega_k^s \\ 0, \mathbf{X} \notin \Omega_k^s \end{cases}. \tag{19}$$

Herein, $A_k = \int_{\Omega_k^s} \mathrm{d}\Omega_k^s$ is the area (volume) of the constructed smoothing domain $\Omega_k^s$.

In accordance with Gaussian divergence theorem, Eq. 16 can be converted to

$$\bar{\boldsymbol{\varepsilon}}_k = \frac{1}{A_k} \int_{\Gamma_k^s} \mathbf{n}_k^s(\mathbf{X}) \mathbf{u}^h(\mathbf{X}) \mathrm{d}\Gamma, \tag{20}$$

where $\Gamma_k^s = \partial \Omega_k^s$, and $\mathbf{n}_k^s$ is a matrix composed of outer normal vector components on the boundary $\Gamma_k^s$, written as (2D)

$$\mathbf{n}_k^s = \begin{bmatrix} n_{kx}^s & 0 \\ 0 & n_{ky}^s \\ n_{ky}^s & n_{kx}^s \end{bmatrix}. \tag{21}$$

Here, $n_{kx}^s$ and $n_{ky}^s$ are the components of unit outward normal vector on $\Gamma_k^s$ to the x- and y-axis, respectively. With the strain smoothing, the classical gradient operation is eliminated in S-FEM.

## 3.2. Formation of ES-FEM

ES-FEM is the most sought-after S-FEM [46]. In this framework, the smoothing domains are constructed based on the edges of elements, as demonstrated in FIG. 2. The initial mesh is generated using three-node linear triangular elements (T3), which simultaneously produces a set of edges. The subsequent task is to identify whether the edge is an inner edge or a boundary edge. Take an internal edge $\Gamma_k^s$ as an example, which is supported by two adjacent elements. Connect the two nodes of the edge to the centroid of its supporting elements, resulting in a quadrilateral smoothing domain, as indicated by the red region in FIG. 2. However, for a boundary edge with only one support element, the smoothing domain is a triangle, as covered by green in FIG. 2.

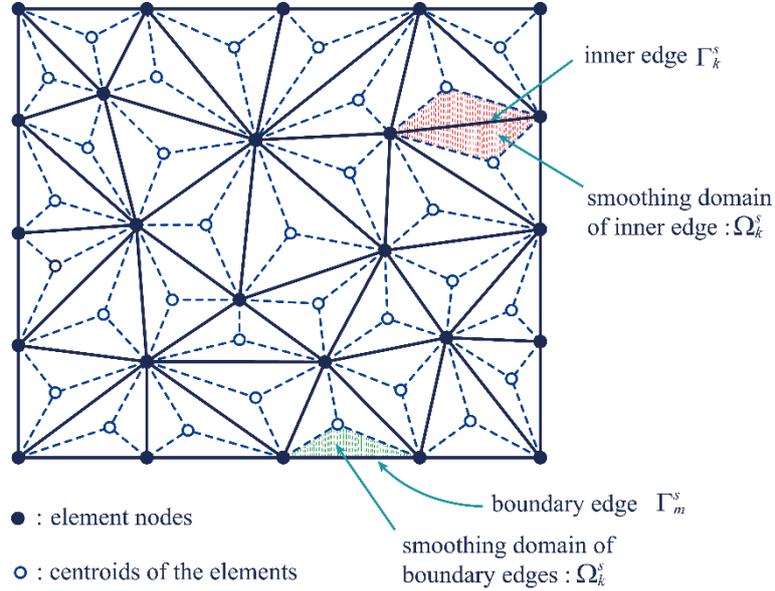

**FIG. 2.** Illustration of smoothing domains in the ES-FEM based on three-node linear triangular elements.

Based on the above scheme, the area of the edge-based smoothing domain can be obtained by

$$A_k^s = \int_{\Omega_k^s} d\Omega = \frac{1}{3}\sum_{i=1}^{N_k} A_i^e, \qquad (22)$$

where $A_i^e$ and $N_k$ are the area and the number of supporting elements, respectively. For the ES-FEM-T3 model currently in use, the smoothing $\bar{\mathbf{B}}_I$ matrix can be straightforwardly given as

$$\bar{\mathbf{B}}_I = \frac{1}{A_k^s}\sum_{i=1}^{N_k}\left[\frac{1}{3}A_i^e \mathbf{B}_i^e\right] \qquad (23)$$

in which, $\mathbf{B}_i^e$ is the standard strain-displacement matrix.

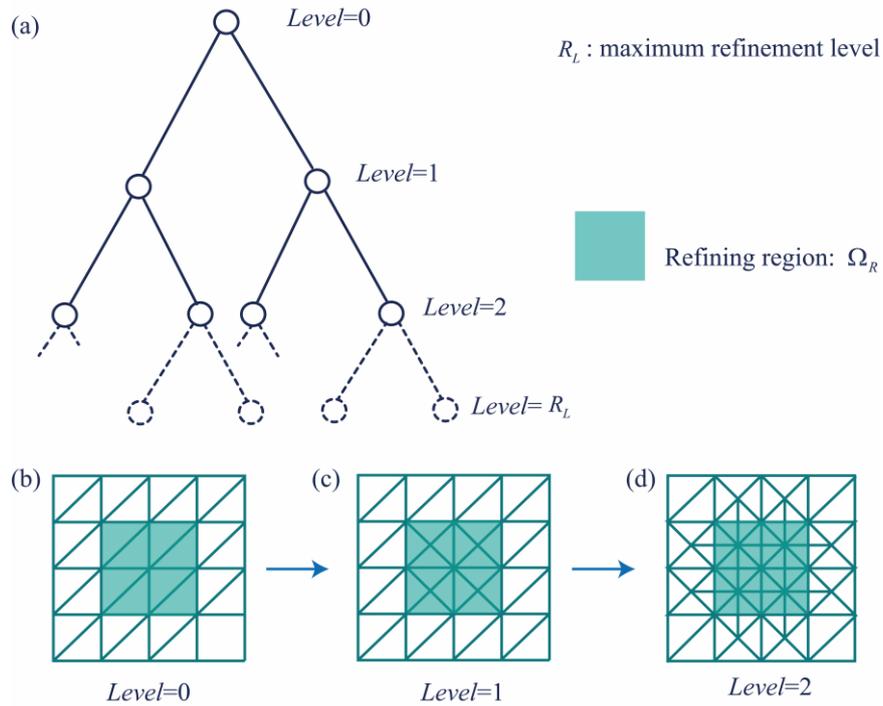

FIG. 3. Illustration of the primary procedures of the proposed adaptive mesh refinement strategy.

3.3 Adaptive mesh scheme

In consideration of the locality of the crack propagation path, phase-field modeling of fracture is a superb application scenario for adaptive mesh algorithm [55, 58, 68]. Thus, an efficient multi-level adaptive mesh strategy developed for the coupling of ES-FEM and PFM was proposed, the primary procedures of which are depicted in FIG. 3. The initial geometry is discretized by triangular elements with the element level ($E_{Level}$) of zero. Before the calculation starts, we specify an element refinement criterion, such as $\phi >= 0.25$, and the maximum refinement level $R_L$. Afterward, all the triangular elements that meet the criteria are bisected into two smaller-scale triangle elements. In case the quality of the refined mesh is poor, an elaborated mesh smoothing scheme based on optimal Delaunay triangulation will be executed [69]. Repeat the above operation until the specified refinement level is reached. Furthermore, the element level is stored according to the binary tree table displayed in FIG. 3(a). And the level of adjacent elements does not differ by more than 1. The corresponding mesh configurations are demonstrated in FIGs. 3(b) to 3(d). It is noteworthy that the refined mesh can be re-coarsened, see our previous work for specific details [57].

After the mesh refinement is accomplished, we map the old nodal data e.g. element connectivity, displacement field and phase-field to the current mesh node. And the edge-based data for ES-FEM like edge connectivity, supporting elements, etc. are also regenerated based on a well-designed program. The outstanding performance of the adaptive mesh presented in the current work will be verified in the Numerical examples section.

## 4. Numerical implementation

### 4.1. Weak form

For ES-FEM, the formation of the weak form is in line with the conventional FEM. Using the standard Galerkin weighted residual method and integration by parts, the weak form of the governing equations given by Eq. 13 can be written as

$$\begin{cases} \int_{\Omega_0} \mathbf{P}:\nabla_0 \delta \mathbf{u} \, dV - \int_{\Omega_0} \bar{\mathbf{b}} \cdot \delta \mathbf{u} \, dV - \int_{\Gamma_N} \bar{\mathbf{t}} \cdot \delta \mathbf{u} \, dA = 0 & \forall \delta \mathbf{u} \in \mathcal{V}_u \\ \int_{\Omega_0} -2(1-\phi) \delta \phi \psi_0 \, dV + \int_{\Omega_0} \left( G_c l_0 \nabla_0 \phi \cdot \nabla_0 \delta \phi + \frac{G_c}{l_0} \phi \delta \phi + \eta \dot{\phi} \delta \phi \right) dV = 0 & \forall \delta \phi \in \mathcal{V}_\phi \end{cases}. \quad (24)$$

Here, $\delta \mathbf{u}$ and $\delta \phi$ are test functions. $\mathcal{V}_u$ and $\mathcal{V}_\phi$ are admissible test spaces for displacement field $\mathbf{u}$ and phase-field $\phi$, which follows

$$\begin{cases} \mathcal{V}_u = \{\delta \mathbf{u} \mid \delta \mathbf{u} = \mathbf{0} \text{ on } \partial \Omega_0 \} \\ \mathcal{V}_\phi = \{\delta \phi \mid \delta \phi = 0 \text{ on } \Gamma_C \} \end{cases}. \quad (25)$$

### 4.2. Linearization of the weak form

Eq. 24 is a set of nonlinear coupled equations, which we first linearize by establishing its Newton-Raphson iteration format, as follows

$$\begin{bmatrix} \frac{\partial \mathbf{R}_u}{\partial \mathbf{u}} & \frac{\partial \mathbf{R}_u}{\partial \phi} \\ \frac{\partial \mathbf{R}_\phi}{\partial \mathbf{u}} & \frac{\partial \mathbf{R}_\phi}{\partial \phi} \end{bmatrix} \begin{bmatrix} \Delta \mathbf{u} \\ \Delta \phi \end{bmatrix} + \begin{bmatrix} \mathbf{R}_u \\ \mathbf{R}_\phi \end{bmatrix} = \begin{bmatrix} \mathbf{0} \\ \mathbf{0} \end{bmatrix}, \quad (26)$$

where the residual vectors $\mathbf{R}_u$ and $\mathbf{R}_\phi$ are defined as

$$\begin{aligned} \mathbf{R}_u &= \int_{\Omega_0} \mathbf{P}:\nabla_0 \delta \mathbf{u} \, dV - \int_{\Omega_0} \bar{\mathbf{b}} \cdot \delta \mathbf{u} \, dV - \int_{\Gamma_N} \bar{\mathbf{t}} \cdot \delta \mathbf{u} \, dA \\ \mathbf{R}_\phi &= \int_{\Omega_0} -2(1-\phi) \delta \phi \psi_0 \, dV + \int_{\Omega_0} \left( G_c l_0 \nabla_0 \phi \cdot \nabla_0 \delta \phi + \frac{G_c}{l_0} \phi \delta \phi + \eta \dot{\phi} \delta \phi \right) dV \end{aligned}. \quad (27)$$

With the non-convexity of energy functional $\prod(\mathbf{u},\phi)$ in mind, the robustness of the monolithic solution of Eq. 26 cannot be guaranteed [32, 35]. A stable staggered algorithm developed by Miehe et al. is therefore employed to decouple the Eq. 26 [23], which is rewritten as

$$\begin{cases} \dfrac{\partial \mathbf{R_u}}{\partial \mathbf{u}} \cdot \Delta \mathbf{u} + \mathbf{R_u} = \mathbf{0} \\ \dfrac{\partial \mathbf{R}_\phi}{\partial \phi} \cdot \Delta \phi + \mathbf{R}_\phi = \mathbf{0} \end{cases} \quad (28)$$

where the directional derivatives $\dfrac{\partial \mathbf{R_u}}{\partial \phi} = \dfrac{\partial \mathbf{R}_\phi}{\partial \mathbf{u}} = \mathbf{0}$. Note that this strategy requires a small load increment to ensure accuracy and robustness.

4.3. Discretization via ES-FEM

As aforementioned, we applied the T3 element for the spatial discretization. The displacement field $\mathbf{u}$, phase field $\phi$ and their gradients can be approximated as

$$\mathbf{u} = \sum_{i=1}^{m} \mathbf{N}_i^{\mathbf{u}} \mathbf{u}_i, \quad \phi = \sum_{i=1}^{m} N_i \phi_i, \quad \boldsymbol{\varepsilon} = \nabla \mathbf{u} = \sum_{i=1}^{m} \mathbf{B}_i^{\mathbf{u}} \mathbf{u}_i, \quad \nabla \phi = \sum_{i=1}^{m} \mathbf{B}_i^{\phi} \phi_i \quad (29)$$

with

$$\mathbf{N}_i^{\mathbf{u}} = \begin{bmatrix} N_i & 0 \\ 0 & N_i \end{bmatrix}, \quad \mathbf{B}_i^{\mathbf{u}} = \begin{bmatrix} N_{i,x} & 0 \\ 0 & N_{i,y} \\ N_{i,y} & N_{i,x} \end{bmatrix}, \quad \mathbf{B}_i^{\phi} = \begin{bmatrix} N_{i,x} \\ N_{i,y} \end{bmatrix} \quad (30)$$

in the framework of normal FEM. Here $N_i$ is the shape function of T3 elements, *m* is the total number of nodes in each element. Whereas in the context of ES-FEM, all the gradient matrices like $\mathbf{B}_i^{\mathbf{u}}$ and $\mathbf{B}_i^{\phi}$ are smoothed in terms of Eq. 23, distinguished by the superscript ($^-$), such as $\bar{\mathbf{B}}_i^{\phi}$.

With respect to the time-dependent dissipative item, a backward Euler difference scheme is utilized, we obtain

$$\dot{\phi} = \frac{\phi_{n+1} - \phi_n}{\Delta t}. \tag{31}$$

Here, $\phi_n$ denotes the phase field at the time $t_n$, and $\Delta t$ is the time step.

Substituting Eqs. 29-31 into Eq. 28, then we have

$$\begin{aligned} \bar{\mathbf{K}}_{uu} \cdot \Delta \mathbf{u} &= -\bar{\mathbf{f}}^{u} \\ \bar{\mathbf{K}}_{\phi\phi} \cdot \Delta \phi &= -\bar{\mathbf{f}}^{\phi} \end{aligned}, \tag{32}$$

in which the smoothing stiffness matrixes $\bar{\mathbf{K}}_{uu}$ and $\bar{\mathbf{K}}_{\phi\phi}$ are written as

$$\begin{aligned} \bar{\mathbf{K}}_{uu} &= \bar{\mathbf{K}}^{mat} + \bar{\mathbf{K}}^{geo} \\ &= \int_{\Omega_k} \bar{\mathbf{B}}_0^T \mathbb{C} \bar{\mathbf{B}}_0 d\Omega + \int_{\Omega_k} \bar{\mathcal{B}}^T \tilde{\tilde{\mathbf{S}}} \bar{\mathcal{B}} d\Omega \\ &= \sum_{k=1}^{N_{sd}} \bar{\mathbf{B}}_0^T \mathbb{C} \bar{\mathbf{B}}_0 A_k^s + \sum_{k=1}^{N_{sd}} \bar{\mathcal{B}}^T \tilde{\tilde{\mathbf{S}}} \bar{\mathcal{B}} A_k^s \end{aligned} \tag{33}$$

and

$$\begin{aligned} \bar{\mathbf{K}}_{\phi\phi} &= \int_{\Omega_k} \left[ \bar{\mathbf{B}}_\phi^T G_c l_0 \bar{\mathbf{B}}_\phi + \bar{N}_\phi (2\psi_0 + \frac{G_c}{l_0}) \bar{N}_\phi^T + \frac{\eta}{\Delta t} \bar{N}_\phi \bar{N}_\phi^T \right] d\Omega \\ &= \sum_{k=1}^{N_{sd}} \left[ \bar{\mathbf{B}}_\phi^T G_c l_0 \bar{\mathbf{B}}_\phi + \bar{N}_\phi (2\psi_0 + \frac{G_c}{l_0}) \bar{N}_\phi^T + \frac{\eta}{\Delta t} \bar{N}_\phi \bar{N}_\phi^T \right] A_k^s \end{aligned}, \tag{34}$$

where the smoothed gradient matrices can be given by

$$\bar{\mathbf{B}}_0 = \begin{bmatrix} \bar{\mathbf{B}}_{I1} \bar{F}_{11} & \bar{\mathbf{B}}_{I1} \bar{F}_{21} \\ \bar{\mathbf{B}}_{I2} \bar{F}_{12} & \bar{\mathbf{B}}_{I2} \bar{F}_{12} \\ \bar{\mathbf{B}}_{I1} \bar{F}_{12} + \bar{\mathbf{B}}_{I2} \bar{F}_{11} & \bar{\mathbf{B}}_{I1} \bar{F}_{22} + \bar{\mathbf{B}}_{I2} \bar{F}_{21} \end{bmatrix} \tag{35}$$

$$\bar{\mathcal{B}} = \begin{bmatrix} \bar{\mathbf{B}}_{I1} & 0 \\ \bar{\mathbf{B}}_{I2} & 0 \\ 0 & \bar{\mathbf{B}}_{I1} \\ 0 & \bar{\mathbf{B}}_{I2} \end{bmatrix} \tag{36}$$

with the smoothed deformation gradient

$$\bar{\mathbf{F}} = \mathbf{I} + \bar{\mathbf{H}} = \mathbf{I} + \begin{bmatrix} \bar{\mathbf{B}}_{I1} \\ \bar{\mathbf{B}}_{I2} \end{bmatrix} \begin{bmatrix} \mathbf{u}_x & \mathbf{u}_y \end{bmatrix}. \tag{37}$$

In the presence of phase-field damage, the tangent modulus (Eq. 8) is multiplied by the degradation function $g(\phi)$, given as

$$\overline{\mathbb{C}} = ((1-\phi)^2 + k)[\mu \overline{J}^{-\beta}[\overline{C}_{il}^{-1}\overline{C}_{kj}^{-1} + \overline{C}_{ik}^{-1}\overline{C}_{jl}^{-1}] + \beta\mu \overline{J}^{-\beta} \overline{C}_{ij}^{-1}\overline{C}_{kl}^{-1}]. \tag{38}$$

Here, $\overline{C}_{ij}$ is the smoothed right Cauchy-Green tensor calculated by

$$\overline{\mathbf{C}} = \overline{\mathbf{H}}^T + \overline{\mathbf{H}} + \overline{\mathbf{H}} \cdot \overline{\mathbf{H}}^T + \mathbf{I}. \tag{39}$$

The matrix $\tilde{\overline{\mathbf{S}}}$ arising in Eq. 33 is defined by

$$\tilde{\overline{\mathbf{S}}} = \begin{bmatrix} \overline{\mathbf{S}} & \mathbf{0} \\ \mathbf{0} & \overline{\mathbf{S}} \end{bmatrix}. \tag{40}$$

Wherein, the PK2 stress combined with the phase-field damage variable $\overline{\mathbf{S}}$ can be written as

$$\overline{\mathbf{S}} = ((1-\phi)^2 + k)\mu(\mathbf{I} - \overline{J}^{-\beta}\overline{\mathbf{C}}^{-1}). \tag{41}$$

Analogous to the construction of the stiffness matrixes, we deduced the following residual vectors in the formation of ES-FEM:

$$\begin{cases} \overline{\mathbf{f}}^{\mathbf{u}} = \sum_{k=1}^{N_{sd}} \overline{\mathbf{B}}_0 \{\overline{\mathbf{S}}\} A_k^s \\ \overline{\mathbf{f}}^{\phi} = \sum_{k=1}^{N_{sd}} \left[ \overline{\mathbf{B}}_\phi^T G_c l_0 \overline{\mathbf{B}}_\phi \phi_{n+1} - \overline{N}_\phi 2\psi_0(1-\phi_{n+1}) + \overline{N}_\phi \frac{G_c}{l_0}\phi_{n+1} + \frac{\eta}{\Delta t}\overline{N}_\phi(\phi_{n+1} - \phi_n) \right] A_k^s \end{cases}, \tag{42}$$

where

$$\{\overline{\mathbf{S}}\} = \begin{bmatrix} \overline{S}_{11} \\ \overline{S}_{22} \\ \overline{S}_{12} \end{bmatrix}. \tag{43}$$

4.4. Irreversibility constraints

Cracks are commonly considered not to heal. Consequently, many methods for enforcing the irreversibility of cracks have been proposed [23, 56, 70]. Among them, a method called activity set is used in the current work [35]. Based on the positive and negative of $\Delta\phi$, the equations related to the calculation of phase-field variable are divided into an active set $\mathcal{I}=\{i|\Delta\phi<0\}$ and its complementary inactive set $\mathcal{I}'$. As a result, we only solve a set of equations reduced according to the inactive set, that is

$$\Delta\phi_{\mathcal{I}'} = -(\mathbf{K}^{\phi\phi})^{-1}_{\mathcal{I}'\mathcal{I}'}(\mathbf{f}^{\phi})_{\mathcal{I}'} \tag{44}$$

with the direct setting $\Delta\phi_{\mathcal{I}} = 0$ in each Newton iteration. The above operation will be executed cyclically until the new set $\{i|\Delta\phi<0\}$ is empty.

### 4.5. Solution procedures

The proposed ES-FEM scheme for phase-field modeling of fracture at large deformation is entirely implemented in MATLAB. Adaptive mesh and time-step algorithms are developed to improve computational efficiency. For clarify, the primary solution procedures are summarized in Algorithm 1.

---

**Algorithm 1**

---

1. Generate initial FEM mesh utilizing T3 elements with $E_{Level} = 0$.

2. Get ES-FEM data based on the initial FEM mesh.

3. Initialize the displacement field $\mathbf{u}_0$ and phase field $\phi_0$ at time $t_0$.

4. Perform a staggered iteration scheme at time step $[t_n, t_{n+1}]$:

    4.1. Initialize the tolerance $tol = 1$.

    4.2. While the tolerance satisfies: $tol \geq 10^{-4}$, run:

        4.2.1. Assemble the smoothed stiffness matrix $\bar{\mathbf{K}}^{\mathbf{uu}}_{i+1}, \bar{\mathbf{K}}^{\phi\phi}_{i+1}$ and residual vectors $\bar{\mathbf{f}}_{\mathbf{u}}, \bar{\mathbf{f}}_{\phi}$ based on the smoothing domains at the iteration step $i+1$.

        4.2.2. Solve $\mathbf{u}^{i+1}_{n+1}$ and $\phi^{i+1}_{n+1}$ based on the Eq. 32 with the input $\mathbf{u}_n$ and $\phi_n$.

        4.2.3. Update the tolerance by

---

$$tol = \max \left\{ \frac{\left\|\bar{\mathbf{f}}^{\mathbf{u}}_{i+1}\right\|}{\left\|\bar{\mathbf{f}}^{\mathbf{u}}_0\right\|}, \frac{\left\|\bar{\mathbf{f}}^{\phi}_{i+1}\right\|}{\left\|\bar{\mathbf{f}}^{\phi}_0\right\|} \right\}. \tag{45}$$

4.3. Output $\mathbf{u}_{n+1}$ and $\phi_{n+1}$ at the current time step.

5. Conduct mesh refinement on the T3 elements whose nodal phase field variable meets $\phi_i \geq 0.25$.

6. Map the old nodal data to the current mesh node.

7. Adaptive adjustment the time step $\Delta t$ and loading step based on the value of $\phi_{n+1} - \phi_n$ (see [35] for detailed algorithms).

8. Advance to the next time step and repeat steps 4-7.

9. Data visualization.

## 5. Numerical examples

In this section, the excellent performance of the proposed method (ES-FEM&APFM) is validated by four representative examples: (i) the double-edge tension specimen with variable length notch; (ii) the fracture of a slab containing a central crack under tension; (iii) the crack propagation of a panel comprising holes; (iv) the crack deflection in hyperelastic materials containing interfaces. Incidentally, the example (i) uses the AT$_2$ model (Eq. 11), and the rest use the AT$_1$ model (Eq. 10). All simulation results are compared with previous experiments and simulations.

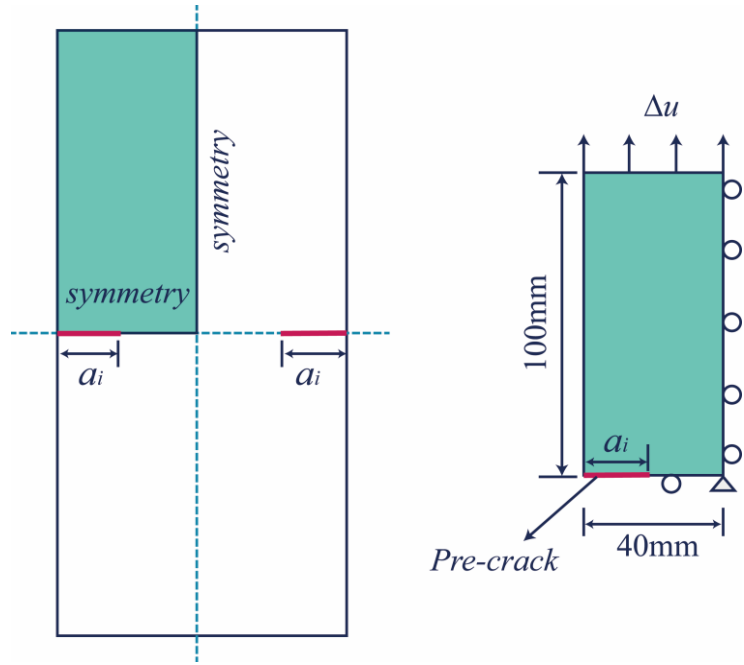

**FIG. 4**. Geometry and boundary conditions of the double-edge notch panel.

5.1. The double-edge tensile specimen with variable length notch

This example originates from the experiment of Hocine et al., which has been regarded as a benchmark [32, 37, 71]. The geometry and boundary conditions are in line with experimental settings, as illustrated in FIG. 4. Utilizing the geometric symmetry, one-quarter of the specimen ( 40 mm×100 mm ) is considered. The initial crack length $a_i$ is varied and the values are taken as $a_i = \{12,16,20,24,28\}$ mm, respectively. For the Neo-Hookean model (Eq. 5), the essential material parameters are set as: $\mu = 0.612$ N/mm$^2$ , $\nu = 0.45$ and $G_c = 7.5$ N/mm . The regularization parameter used in this simulation is set to $l_0 = 1$ mm , correspondingly, the effective size of the element in the vicinity of the crack path is determined by $h_f \approx l_0/8$ . Benefiting from the adaptive mesh algorithm, the total amount of elements is only 3092 in the initial stage, and it eventually reaches 18532 with crack growth. Due to numerical stability concerns, a non-zero viscosity coefficient $\eta = 1 \times 10^{-3}$ and adaptive time step and loading step algorithms are introduced, which can draw on the work of Loew et al. [35]. The crack evolution patterns and the corresponding adaptive mesh with the pre-crack length $a_i = 16$ mm at five different loading steps are shown in FIG. 5. To visualize the crack opening at large deformation, the level set of phase field variable that satisfies $\phi > 0.8$ is removed in the current configuration. From the crack initiation (FIG. 5(b)) to the complete break (FIG. 5(e)), our results are a strong resemblance to those in the literature (a complete movie is available in the Supplemental Materials). A noteworthy fact is that the proposed adaptive ES-FEM can considerably improve computational efficiency, which is foreseeable in view of FIGs. 5(f)-(g). In the present test, it's about 20 times faster than the standard ES-FEM. Besides, we also recorded the load-displacement curves with varying pre-crack length throughout the loading stage, as depicted by the solid lines in FIG. 6. Evidently, the simulation results match well with the experimental measurements extracted from the literature [71] (the dash lines in FIG. 6).

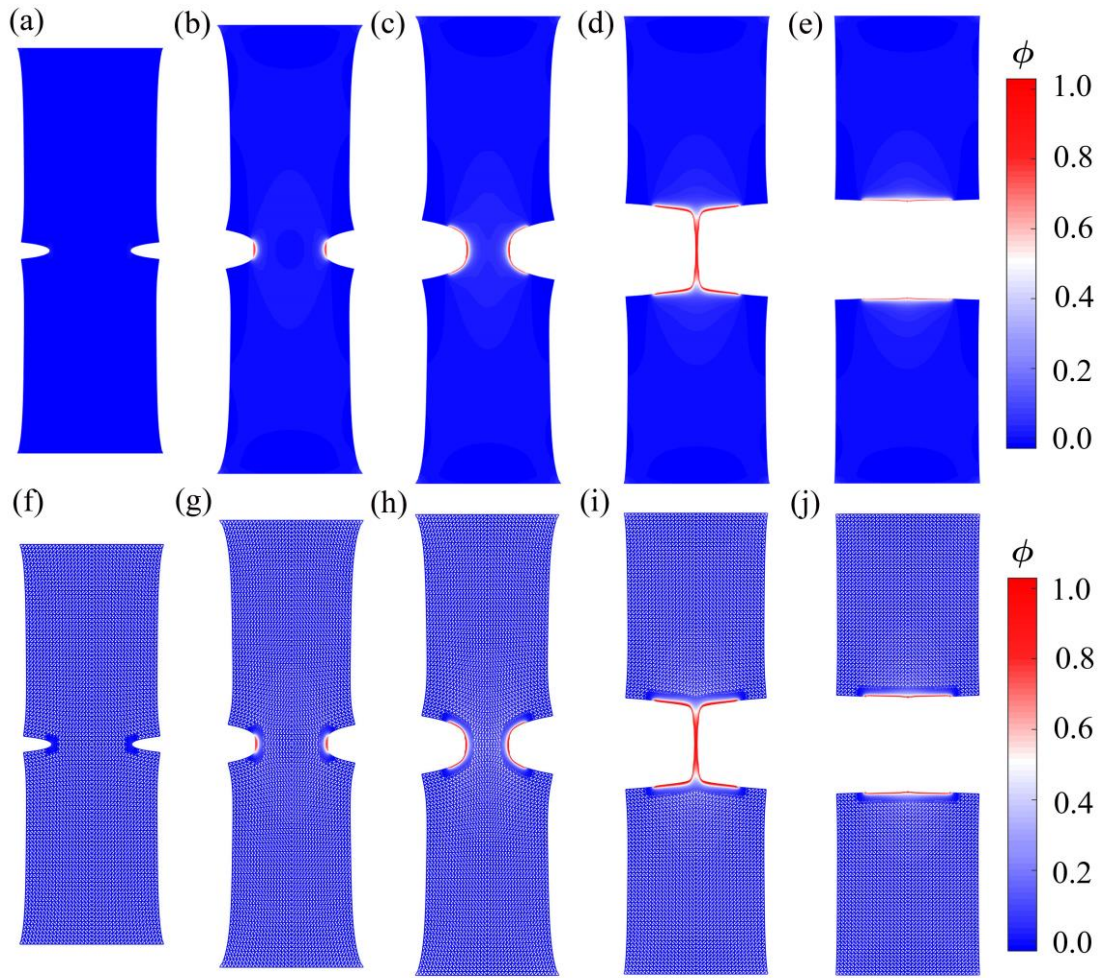

**FIG. 5.** Crack evolution patterns and the corresponding adaptive meshes of the double-edge tensile specimen with the pre-crack length $a_i = 16$ mm at loading displacements $\Delta u = 25.0000$ mm, $\Delta u = 56.0000$ mm, $\Delta u = 58.2610$ mm, $\Delta u = 58.2614$ mm, $\Delta u = 58.2616$ mm, respectively.

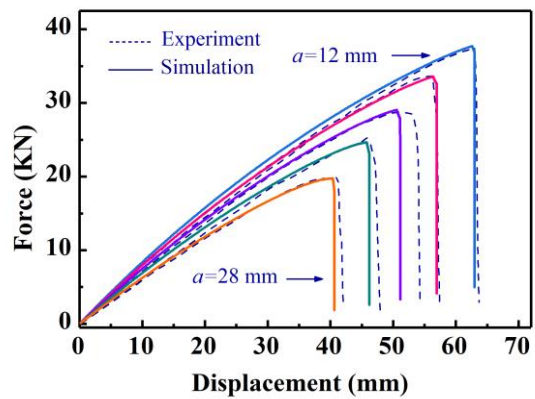

**FIG. 6.** Comparison of simulation results and experiments of force-displacement curves for the double-edge tensile specimen with five different notch length: $a_i = 12, 16, 20, 24, 28$ mm.

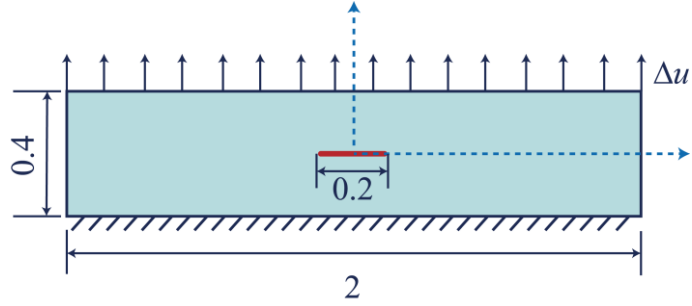

**FIG. 7.** Illustration of the initial geometry and applied boundary conditions for the slab containing a central crack tension test (sizes in mm).

5.2. The fracture of a slab containing a central crack under tension

In this test, the initial geometry and boundary conditions are established in accordance with the benchmark [32, 54], as demonstrated in FIG. 7. Likewise, only a quarter of the specimens were involved in the calculation allow for the symmetry. The loading controlled by the vertical displacement is adaptively modified with the step increment of the phase field. For the sake of comparison with the results of Kumar et al., the constitutive model (Eq. 5) is rewritten as

$$\psi_0(\mathbf{F}) = \frac{\mu}{2}[\text{tr}[\mathbf{F}^T\mathbf{F}] - 3] + \frac{\mu^2}{\Lambda}[(J^{-\Lambda/\mu} - 1)], \tag{46}$$

where the material parameters are set as $\mu = 5$ N/mm$^2$ and $\Lambda = 7.5$ N/mm$^2$. Moreover, the critical fracture energy $G_c$ is 3 N/mm and the regularization parameter is taken as $l_0 = 0.01$ mm. Like the previous example, we also implemented the adaptive ES-FEM algorithm in which the effective element size is approximate $l_0/10$. The morphology of the cracks and their corresponding adaptive meshes from sprouting to complete fracture at three different deformation states ($\Delta u = 0.2600$ mm, $\Delta u = 0.4680$ mm, $\Delta u = 0.4736$ mm) is presented in FIG. 8 (An animation of the complete crack evolution for this example is provided in the supplemental materials). Intuitively, the crack patterns are almost identical to those in the literature. For comprehensive consideration, the resulting force-displacement curve (labeled by ES-FEM&APFM) is plotted and compared to the literature [2]. As indicated in FIG. 9, their matching is clearly satisfactory. It is worth mentioning that the standard FEM converges slowly and precariously in the same mesh situation. However, the proposed method (ES-FEM&APFM) does not encounter this problem.

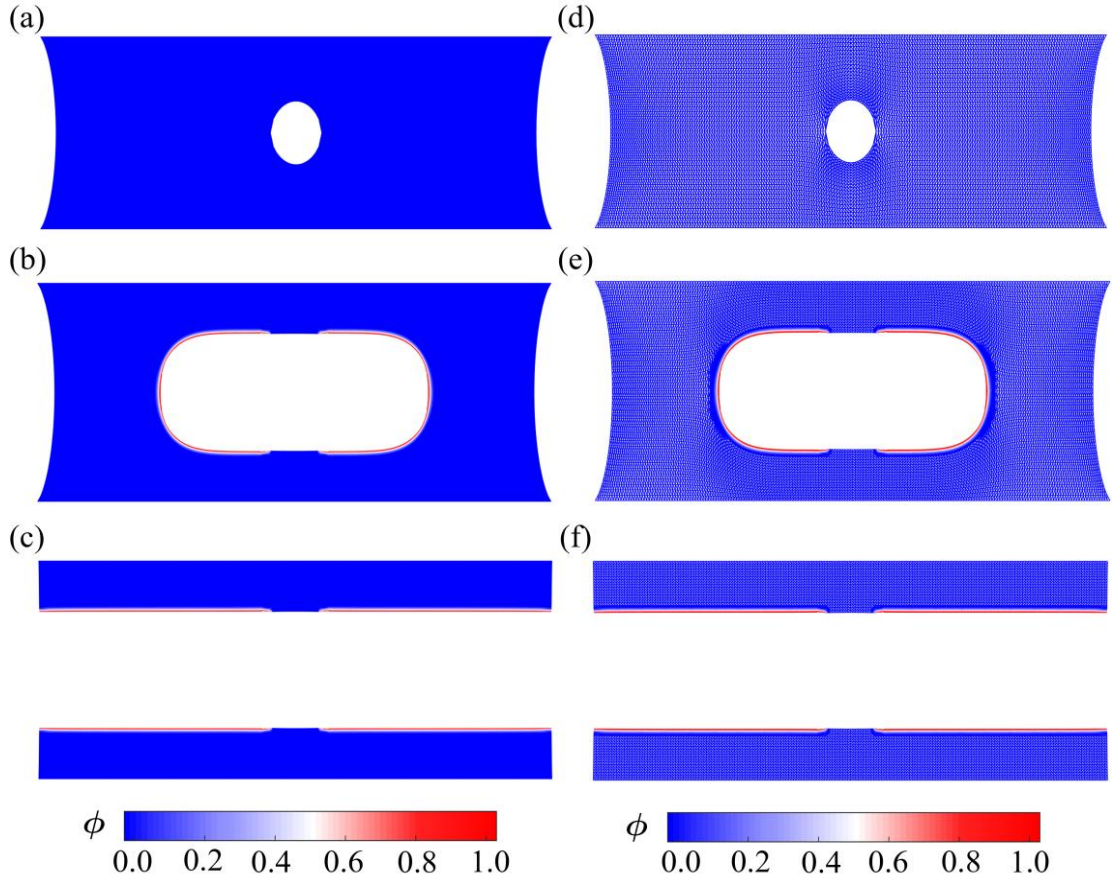

**FIG. 8.** The morphology of the cracks and their corresponding adaptive meshes for the slab containing a central crack tension test at loading displacement $\Delta u = 0.2600$ mm, $\Delta u = 0.4680$ mm, $\Delta u = 0.4736$ mm, respectively.

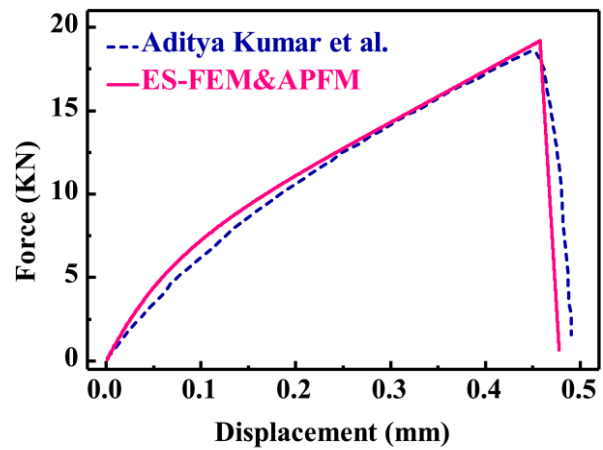

**FIG. 9.** Load-displacement curves of the slab containing a central crack tension test.

## 5.3. The crack propagation of a panel comprising holes

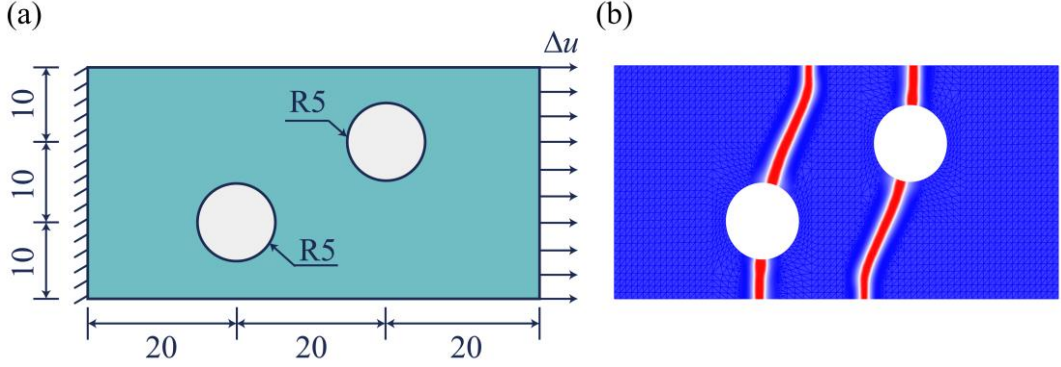

FIG. 10. (a) Illustration of the initial geometry and boundary conditions for the tension test of the panel comprising holes (sizes in mm). (b) Crack patterns in the context of the initial grid.

This example involves the propagation of the curvilineal cracks deflected by two eccentric holes at large deformation [72]. The initial geometry and the applied boundary conditions are illustrated in FIG. 10(a). And the final crack patterns in the context of initial configuration are presented in FIG. 10(b). An adaptive horizontal displacement loading is imposed on the right edge while the left edge is constrained. We assume the material parameters are $\mu = 0.28$ N/mm$^2$, $\nu = 0.45$, $G_c = 1.4$ N/mm, $\eta = 1\times 10^{-3}$ and $l_0 = 0.5$ mm. The size of the fine-scale element in the adaptive mesh approximately satisfies $h_f = l_0/5$. Note that we did not set an explicit pre-crack in this example. Along with the applied constant displacement incremental loading accumulates, the cracks initiate immediately after the stress around the holes exceeds a threshold $\sigma_c$, as shown in FIG. 11(a). Thereafter, a self-regulating slow loading is applied until completely fracture. Snapshots of crack patterns in several deformed states are displayed in FIGs. 11(b)-(e), and the corresponding animation of the crack evolution is provided in the supplemental materials. Furthermore, FIGs. 11(f)-(j) present the evolution of the corresponding adaptive mesh, clearly demonstrating the effectiveness of the proposed algorithm. The entire crack propagation patterns for this example are similar to that of Miehe et al. [72], such as the final three fragments. However, quantitative comparisons have not been conducted in view of the lack of available force-displacement curves.

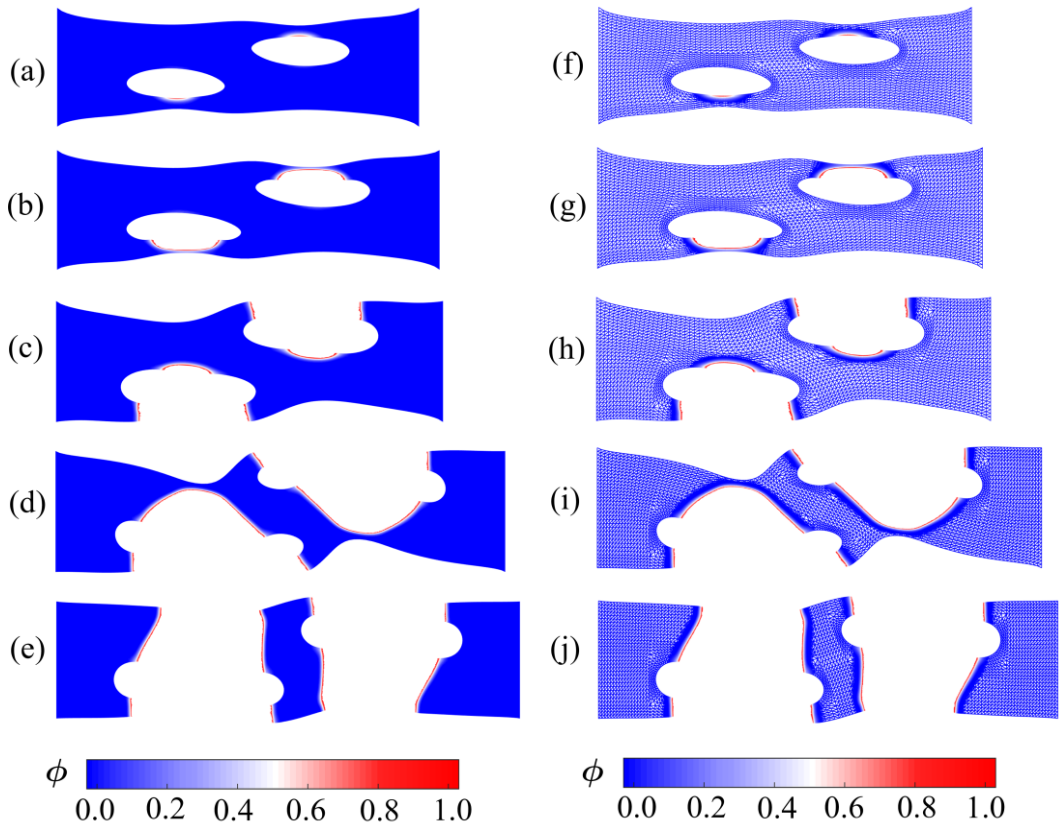

FIG. 11. Snapshots of crack patterns and the corresponding adaptive meshes for a panel comprising holes tension test at loading $\Delta u = 17.5000$ mm, $\Delta u = 17.6050$ mm, $\Delta u = 17.6057$ mm, $\Delta u = 27.1451$ mm, $\Delta u = 31.1287$ mm.

5.4. The crack deflection in hyperelastic materials containing interfaces

This test is designed to investigate the deflection effect of weak interfaces on crack propagation in hyperelastic materials such as rubber and hydrogel [73, 74]. At first, we consider a rectangular strip with a width of 24 mm and a height of 120 mm, as illustrated in FIG. 12(a). A horizontal notch of 12 mm in length is cut from the left edge at the middle height of the specimen. A symmetrical displacement loading is imposed on the upper and lower edges. At the horizontal symmetry center, a weak interface of 0.8 mm in width is created by specifying $G_c^{bulk}/G_c^{interface} \geq 10$ based on local refinement mesh. Herein, $G_c^{bulk}$ and $G_c^{interface}$ denote the critical fracture energy of bulk and interface, respectively. The essential material parameters for this example are set to $\mu = 0.035$ N/mm$^2$, $\nu = 0.45$, $G_c^{bulk} = 0.034$ N/mm, $G_c^{interface} = 0.0017$ N/mm, $\eta = 1 \times 10^{-3}$ and $l_0 = 0.2$ mm. FIG. 12(b) presents the crack patterns at three different loading stages (a complete crack evolution movie can be found in supplemental materials). As observed, instead of direct penetration, a straight crack deflects at the weak interface. Intriguingly, this result agrees well with the experimental photographs shown in FIG. 12(c) [74]. Some simulation results not presented here reveal that the crack penetration vs deflection is largely relying on the ratio of $G_c^{bulk}/G_c^{interface}$, nevertheless, we do not intend to conduct an in-depth discussion in the current work.

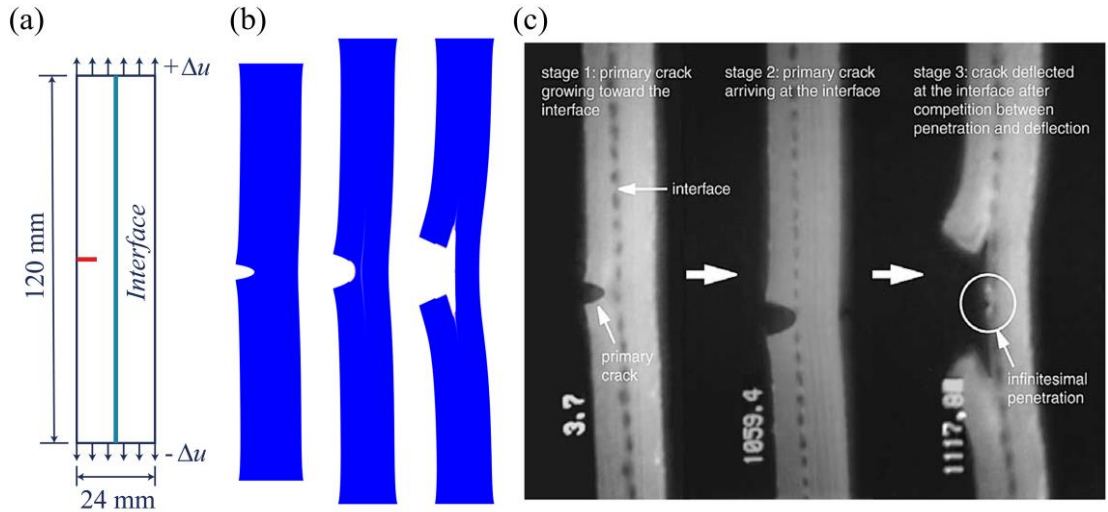

**FIG. 12.** Crack deflection in hyperelastic materials containing interfaces. (a) Initial geometry and boundary conditions. (b) Crack patterns at three different deformation stages. (c) Experimental snapshots of bi-layered rubber interface fracture.

## 6. Conclusion

In summary, phase-field modeling of fractures at large deformations was first formatted in the ES-FEM framework. The nature of the phase-field approach endues it inimitable advantages in simulating fractures. ES-FEM is developed by combining FEM and meshless ideas with high accuracy and insensitivity to element distortion. The current work combines PFM with ES-FEM, which initiates a novel approach for modeling the fracture at large deformations. However, PFM generally requires fine meshing to correctly identify the crack trajectory, and the bandwidth of the stiffness matrix for the ES-FEM is larger than that of the conventional FEM. Therefore, albeit the combination of the two methods manifests high precision, faster convergence rate and better robustness than existing approaches, it was testified as computationally demanding. For this reason, a multi-level adaptive mesh scheme designed for the coupling of ES-FEM and PFM was presented. In addition, we also outlined the specific numerical implementation of the proposed method (ES-FEM&APFM), the effectiveness of which is verified by several representative numerical examples. Particularly, the experiment of weak interface frustration crack propagation in a rubber-like solid was first reproduced by our approach. In the next work, we will conduct a thorough investigation on the effects of weak interfaces on the competition of crack penetration vs deflection at large deformations.

**ACKNOWLEDGMENTS**

The authors would like to thank Prof. Shan Tang (Dalian University of Technology) for his valuable suggestions on the numerical algorithm of PFM at large deformations. This work is supported the National Natural Science Foundation of China (51633009, and 51790500).